\documentclass[11pt]{amsart}
\usepackage[makeroom]{cancel}
\usepackage{color}
\usepackage{amsmath,amssymb,amsthm}

\newtheorem{theorem}{Theorem}[section]
\newtheorem{lemma}[theorem]{Lemma}

\newtheorem{assumption}[theorem]{Assumption}
\newtheorem{remark}[theorem]{Remark}

\newcommand{\A}{\mathcal{A}}
\newcommand{\D}{\mathcal{D}}
\newcommand{\F}{\mathcal{F}}
\newcommand{\cH}{\mathcal{H}}

\newcommand{\R}{\mathcal{R}}
\newcommand{\Ahalf}{\mathcal{A}^{1/2}} 

\begin{document}
\title[MGT equation with memory]{Moore-Gibson-Thompson equation with memory, part I: exponential decay of energy}

\author{Irena Lasiecka}
\address{%
Department of Mathematical Sciences\\
University of Memphis, Memphis, TN, 38152}
\email{lasiecka@memphis.edu}
\author{Xiaojun Wang*}
\address{%
Department of Mathematical Sciences\\
University of Memphis, Memphis, TN, 38152}
\email{xwang13@memphis.edu}

%
%
\maketitle
\begin{abstract}
We are interested in the Moore-Gibson-Thompson(MGT) equation with memory \begin{equation}\nonumber
\tau u_{ttt}+ \alpha u_{tt}+c^2\A u+b\A u_t -\int_0^tg(t-s)\A w(s)ds=0. \end{equation}
We first classify the memory into three types. Then we study how a memory term creates damping mechanism and how the memory causes energy decay.\end{abstract}

\tableofcontents

\section{Introduction}

We are interested in the Moore-Gibson-Thompson(MGT) equation with memory
\begin{equation}\label{mainmgt}
\tau u_{ttt}+ \alpha u_{tt}+c^2\A u+b\A u_t -\int_0^tg(t-s)\A w(s)ds=0,
\end{equation}
where $\tau, \alpha, c^2, b$ are physical parameters and $\A$ is a positive self-adjoint operator on a Hilbert space $H$. The convolution term $\int_0^tg(t-s)\A w(s)ds$ reflects the {\it memory effects} of materials due to viscoelasticity. Here the convolution kernel $g$ satisfies proper conditions exhibiting  ``memory character", which will be explained later. The form of $w$
classifies the memory into the following three types
\begin{eqnarray}\label{mgtm}
&&w(s)=u(s), \hspace{1.05in}\text{\bf(Type 1)}\\
&&w(s)=u_t(s),\hspace{1in} \text{\bf (Type 2)}\\
&&w(s)=\lambda u(s)+u_t(s). \hspace{.5in}\text{\bf (Type 3)}
\end{eqnarray}

This is the first part of an ongoing project, in which a {\it third order (in time)} equation with memory is investigated.  More specifically, we study how (and whether) a memory term creates damping mechanism and how (and whether) it causes energy decay.

\subsection{Background}
Moore-Gibson-Thompson (MGT) equation arises from modeling high frequency ultrasound waves. Without memory, the linearized MGT equation reads
\begin{equation}\label{mainmgt0}
\tau u_{ttt}+ \alpha u_{tt}+c^2\A u+b\A u_t =0.
\end{equation}
Certainly this equation is in abstract form and it has a simple prototype where $\A=-\Delta$.
In \cite{KLM11}, the well-posedness of (\ref{mainmgt0}) and the uniform decay of its energy is studied under proper functional setting and initial-boundary conditions. We shall revisit the results of \cite{KLM11} in Section 2.
Spectral analysis for this model has been carried out in \cite{MMT12}, which confirms the validity and sharpness of the results in \cite{KLM11}.

A linear MGT equation is the prelude to nonlinear ones. The classical nonlinear acoustics models include the Kuznetsov equation, the Westervelt equation and the KZK equation, etc. This research field is highly active due to a wide range of applications such as the medical and industrial use of high intensity ultrasound in lithotripsy, thermotherapy, ultrasound cleaning, etc.  There have been quite a few works in this aspect, more from engineering viewpoint; see \cite{Jor09, KLM11, MG60, Kuz71, Cri79}. Rigorous study of the small global solutions and their decay estimates has been conducted in \cite{KLP12}. A thorough study of the linearized models is a good starting point for better understanding the nonlinear models. Actually, the work \cite{KLM11} has shown that, even in the linear case, rich dynamics appear.

Here we investigate the linear MGT system when memory effects get involved. For reasons that will be clear soon, we start with introducing {\it wave equation}. 

It is well known that wave system conserves mechanical energy. Consider the wave equation
$$u_{tt}-\Delta u=0.$$ 
Multiplying it by $u_t$, which is called a ``multiplier", we have
$${d\over dt}{1\over 2}(||u_t||^2+||\nabla u||^2)=0\Rightarrow E(t)\equiv{1\over 2}(||u_t||^2+||\nabla u||^2)=E(0).$$
Here $E(t)$ represents the mechanical energy, summation of kinetic and potential energy. We see the mechanical energy does not change when a wave evolves. On the other hand, we normally see mechanical energy decays in a physical system due to certain damping mechanism. A different viewpoint is, in order to force the energy decay, we have to inject damping mechanisms into the system. There are different ways to implement this idea. One way\footnote{Adding structure damping $-\Delta u_t$ to the wave equation is another way to obtain exponential decay of the energy. These topics are of great importance in real applications, hence forming an active research area \cite{Adh00}.
} is to add viscous damping, also called frictional damping, into the wave equation. Taking into account the friction, we end up with equation
$$u_{tt}-\Delta u+u_t=0.$$ 
The same multiplier $u_t$ gives energy estimate $${d\over dt}{1\over 2}(||u_t||^2+||\nabla u||^2)=-||u_t||^2\leq 0\Rightarrow E(t)\leq E(0).$$ While the energy functional $E(t)$ keeps the same form ${1\over 2}(||u_t||^2+||\nabla u||^2)$, the multiplier $u_t$ generates a pure negative term $-||u_t||^2$ on the right, causing damping in the system. A quick thought would be, if $||u_t||^2$ controls the full energy, say $||u_t||^2\geq CE(t)$, then we have ${d\over dt}E(t)\leq -CE(t)$ and Gronwall inequality immediately gives the exponential decay of energy, $E(t)\leq C_1e^{-Ct}$. It is not the case here, however. 

Multiplier $u_t$ only gives us partial information: the mechanical energy is decreasing (at least when $u_t\neq 0$). Without further investigation, It is not clear how fast it decreases, or whether it decreases exponentially. Fortunately, the energy here does decay exponentially. Indeed, it can be shown, by standard Lyapunov function method (essentially with the help of another multiplier $u$), that the energy decays exponentially \cite{Che79}. 

What interests us here is a different kind of damping mechanism, namely the energy dissipation caused by viscoelasticity, which forces the appearance of memory term in the system. Viscoelasticity is  the property of materials that exhibit both viscous and elastic characteristics when undergoing deformation. It usually appears in fluids with complex microstructure, such as polymers. One encounters viscoelastic materials in biological science, materials sciences as well as in many industrial processes, e.g., in the chemical, food, and oil industries. The phenomena and mathematical models for such materials are more varied and complex than that of pure elastic materials or that of pure Newtonian fluids, see\cite{Ren00}.

Our work is motivated by a series of papers on the wave equation (or its abstract version), with memory 
\begin{equation}\label{wave_memory}
u_{tt}-\Delta u+\int_0^tg(t-s)\Delta u(s)ds=0,
\end{equation}
where $\Delta$ is the Laplacian operator defined on certain functional space. The convolution term $\int_0^tg(t-s)\Delta u(s)ds$ represents memory: the integral itself suggests the nonlocality in time; the system at present moment ``remember" information in the past history. Generally $g(t)$ is assumed to be decreasing, suggesting the ``nearer past" has stronger influence on the system than the ``further past" does. This is a simplified one dimensional wave equation with memory. Like all the physically meaningful partial differential equations, it is derived based on the balance laws of momentum and/or conservation of mass, together with a constitutive equation relating stress to strain (elasticity) or stress to rate of strain (fluids). For more information on the derivation of this 1-D system, we refer to \cite{HNR87}.

As for {\it wave equation}, there are hundreds of papers studying either the damping effects caused by the memory terms or the interaction between different damping mechanisms. See \cite{ACS08, AC09, Ala10, Ala12, CCLW15, CCM08, FP02, HW10, LMM13, LW14, Mes08, MM12, RS01} to name a few.
While a memory term does generate damping in a wave system, the nature of this mechanism seems to be subtle. In \cite{FP02},  the authors found that a damped wave equation, which is already exponentially stable, loses this property after being added a memory that is not properly calibrated. This suggests, {\it the more is not necessarily the better.} It is important to understand the potential benefits as well as drawbacks of adding memory to a system in applications. 

It is natural to ask questions such as, what kind of term creates damping effect in a system? is it possible that adding a damping term ``hurt" the energy, in any sense? The answer hides in the energy equality in a wave-memory system. Roughly speaking, the memory term creates damping by weakening the total energy in the first place: it borrows energy to initialize the deformation and pays back in the long run. We exploit this point in another paper \cite{LW15}.

To our best knowledge, the present paper is the first work to investigate a temporally {\it third order} system that interacts with memory. As we will see in Section 2, the MGT equation is separated into different regimes based on a critical parameter $\gamma=\alpha-{c^2\tau\over  b}$. In the non-critical regime($\gamma>0$), which will be explained presently, a substantial portion of frictional damping exists. Hence it is not surprising that the energy in non-critical case intends to decay more often than that in critical case. 

A big part of the work is on how the memory term influences the MGT system in the {\it non-critical case}, namely for $\gamma\neq 0$. This is not particularly challenging since the system alone contains strong damping mechanism. But this study brings up more interesting questions such as how different damping mechanisms ``cooperate" with each other. This is never a trivial issue, see \cite{CO03,LW15}. Moreover, this step provides intuitions when it comes to the {\it critical case} $\gamma=0$, where the system itself does not provide enough damping. In the {\it critical case}, we show that ``the right mixture" of memory can make the MGT system exponentially stable\footnote{It is currently an open question whether the first type of memory is sufficient to make the energy decay.}. Based on the results in wave system, this result is surprisingly neat, see Theorem \ref{theorem32}.

The following notations are frequenctly used in the work.
\begin{itemize}
\item $\A$: a positive selfadjoint operator defined on a real Hilbert space $H$ with domain $\D(\A)$. We assume $\lambda_0>0$ exists such that $||u||\leq \lambda_0||\Ahalf u||$ for all $u\in H$.
\item$|| u|| \triangleq ||u||_H $: the norm on $H$.
\item $\cH\triangleq\D(\Ahalf)\times\D(\Ahalf)\times H$.
\item
$g\circ u \triangleq \int_0^tg(t-s) ||u(t)-u(s)||^2 ds$: an important quantity in the energy estimates.
\item
$g* u (t)\triangleq\int_0^t g(t-s) u(s) ds $: the convolution integral.
\item 
$G(t)\triangleq\int_0^t g(s)ds$: the cumulative strength of viscoelasticity.
\item 
$F_0(t)\triangleq ||u_{tt}||^2 + ||\A^{1/2} u_t||^2 + ||\A^{1/2} u||^2$: the standard energy defined for MGT.
\item
$F_1(t)\triangleq F_0(t)  - g' \circ \A^{1/2} u$: the standard energy for MGT with memory of type 1.
\item 
$F_{2} (t)\triangleq F_0(t)  + g \circ \A^{1/2} u_t$: the standard energy for MGT with memory of type 2.
\item
$F_3(t)\triangleq F_0(t)+g\circ \Ahalf w$: the standard energy for MGT with memory of type 3, with $w=\lambda u+u_t$.
\item
$ F_3^{cr}(t)\triangleq ||\A^{1/2} w||^2 + ||w_t||^2 + g \circ  \A^{1/2} w$: energy for type 3 memory, in critical regime.
\item $\alpha, \tau, c^2, b$: parameters in MGT equation.
\item $\gamma=\alpha-{c^2b\over \tau}$: a critical paremeter in defining concepts of non-critical and critical regimes.
\item $c_0, c_1$, ..., etc: constants related to memory kernel $g$.
\item $C, C_0, C_1$,...,etc: generic constants related to energy estimate.
\end{itemize}

\subsection{Main Results} 


We make the standing assumption
\begin{assumption}\label{assumption0} 
Assume 
\begin{itemize}
\item 1. $g(\cdot) \in C^1(0,\infty)\bigcap C[0,\infty)$, and  $g(t)\geq 0, g'(t) \leq 0$ for $ t\geq 0$. 
\item 2. $g''(t)\geq 0$ for $ t\geq 0$. 
\item 3. Positive constants $c_0, c_1$ exist such that $g'\leq -c_0g$.
\end{itemize}
\end{assumption}
\begin{remark}
Among the above list, item 1 is the standing assumption on memory kernel. Item 2 says the kernel must be convex. This restriction can actually be removed; but to get energy bound, we need to impose condition on g(0).  see \cite{LW15}. Item 3 implies the exponential decay of the kernel. In \cite{LW15}, we deal with problem that has memory kernel with general decay rate.
\end{remark}

Next we present the main results of energy decay in case of different memory type. \\

{\bf *Memory of type 1:}\\

We consider
$$\tau u_{ttt}+ \alpha u_{tt}+c^2\A u+b\A u_t -\int_0^tg(t-s)\A u(s)ds=0.$$

We make the following assumption on the memory kernel.
\begin{assumption}\label{assumption1} 
Assume 
\begin{itemize}
\item 1. $\gamma=\alpha-{c^2\tau\over b}>0$.
\item 2. $G(+\infty)<c^2$.
\end{itemize}
\end{assumption}

\begin{theorem}\label{theorem1}
Under Assumptions \ref{assumption0} and \ref{assumption1}, the unique weak solution $u$ exists on $\cH$ and the standard energy $F_1(t) $ decays exponentially, i.e.,  there exist positive constants $ C, \omega > 0 $ such that 
$$ F_1(t) \leq CF_1(0) e^{-\omega t },  t > 0.$$
\end{theorem}


{\bf *Memory of type 2:}\\

We consider
$$\tau u_{ttt}+ \alpha u_{tt}+c^2\A u+b\A u_t -\int_0^tg(t-s)\A u_t(s)ds=0.$$
%

\begin{assumption}\label{assumption2} 
Assume 
\begin{itemize}
\item 1. $\gamma=\alpha-{c^2\tau\over b}>0$.
\item 2. There exist positive numbers $\theta$ and $k \in (\frac{c^2}{b}, \frac{\alpha}{\tau} )$ such that  ${k\over \theta}<c_0$ and 
$$G(+{\infty}) \leq \min\{{2(bk-c^2)\over 2+\theta}, b-{c^2\over k}\}.$$
\end{itemize}
\end{assumption}

\begin{theorem}\label{theorem2}
Under Assumptions \ref{assumption0} and \ref{assumption2}, the unique weak solution $u$ exists on $\cH$ and the energy satisfies $$ F_2(t) \leq CF_2(0) e^{-\omega t },  t > 0.$$
\end{theorem}

{\bf *Memory of type 3:}\\


We consider MGT with a mixed memory in form of
$$\tau u_{ttt}+ \alpha u_{tt}+c^2\A u+b\A u_t -g*\A w=0, $$
where $w=\lambda u(t)+u_t(t).$

For the {\bf non-critical case}, $\gamma=\alpha-{c^2\tau\over b}>0$. Define the standard energy functional $$F_3(t)=||\Ahalf u_t||^2+||\Ahalf u||^2+||u_{tt}||^2+g\circ \Ahalf w.$$

\begin{assumption}\label{assumption31} 
Assume 
\begin{itemize}
\item 1. $\gamma=\alpha-{c^2\tau\over b}>0$.
\item 2. $\lambda\in (\frac{c^2}{b}, \frac{\alpha}{\tau})$ such that $G(+\infty)<{c^2\over \lambda}$.
\end{itemize}
\end{assumption}

\begin{theorem}\label{theorem31}
Accept the Assumptions \ref{assumption0} and \ref{assumption31}, then the unique weak solution $u$ exists on $\cH$ and the energy decays exponentially, i.e.
$$F_3(t) \leq CF_3(0) e^{-\omega t },  t > 0.$$
\end{theorem}

For the {\bf critical case} where $\gamma=0$, the energy corresponding to the above system becomes 
$$ F_3^{cr}(t)  \triangleq ||\A^{1/2} w||^2 + ||w_t||^2 + g \circ  \A^{1/2} w.$$

\begin{assumption}\label{assumption32} 
Assume 
\begin{itemize}
\item 1. $\gamma=\alpha-{c^2\tau\over b}=0$.
\item 2. $\lambda={\alpha\over\tau}, G(+\infty)<{c^2\over \lambda}$.
\end{itemize}
\end{assumption}

\begin{theorem}\label{theorem32}
Accept the Assumptions \ref{assumption0} and \ref{assumption32}, then the unique weak solution $u$ exists on $\cH$ and the energy $F_3^{cr}(t)$ decays exponentially, namely we have the exponential decay of $||w_t||$ and  $||\Ahalf w||$.
\end{theorem}

\begin{remark}
The existence part of each proof shall be omitted; we focus on the study of decay rate. In fact the existence can be easily carried out by the standard Galerkin method. The main ingredient of Galerkin method is certainly the global a priori bounds of the solutions, which can be found in the stability results.
\end{remark}

\begin{remark}
In this work, we only consider exponential decay, which demands the exponential decay of the memory kernel, namely $g'\leq -c_0 g$. In paper \cite{LW15}, based on the idea in \cite{LT93}, we will deal with the more general case $g'+H(g)\leq 0$, where $H(\cdot)$ is a convex increasing function with $H(0)=0$.
\end{remark}

The paper is organized as follows:  In Section 2 we revisit the linearized Moore-Gibson-Thompson equation. The main result from \cite{KLM11} is stated and a slight different proof is presented. In Sections 3, 4 and 5, we give the proofs of the main theorems stated above.

\section{Revisit the MGT equation}

The MGT model is third order in time. The relaxation parameter $\tau$ accounts for finite speed of acoustic waves, addressing the paradox of infinite speed of propagation occurring in modeling acoustic waves \cite{KLM11}. 
As a consequence of the presence of $\tau$, the original parabolic like model($\tau =0$) becomes hyperbolic($\tau>0$); the issues of well-posedness, stabilization and long time behaviors become more involved due to the presence of infinitely many unstable eigenvalues. With $\tau > 0$, roughly speaking when diffusion parameter $b=0$, (\ref{mainmgt0}) becomes ill-posed;  on the other hand with $b > 0$ the system generates a strongly continuous semigroup. Moreover, the corresponding semigroup is exponentially stable on $\cH$ when $\gamma=\alpha - c^2\tau b^{-1} > 0$. When $\gamma =0$ the semigroup is conservative.

More specifically, consider the MGT equation (\ref{mainmgt0}),
$$\tau u_{ttt}+ \alpha u_{tt}+c^2\A u+b\A u_t=0,$$ with initial conditions
$$u(0)=u_0, u_t(0)=u_1, u_{tt}(0)=u_2.$$
Here $\A$ is a self-adjoint positive operator on a Hilbert space $H$ with a dense domain $\D(\A)\subset H$. The initial data $$(u_0, u_1, u_2)\in \cH\equiv\D(\Ahalf)\times\D(\Ahalf)\times H.$$
Here the parameters $\tau>0, b>0, \alpha\in\R$ represent the relaxation, the diffusivity and the friction. We state the following results from \cite{KLM11}:

\begin{theorem}[Kaltenbacher et.al. 2011, \cite{KLM11}]\label{KLM11_1}

Let $b>0, \alpha\in \R$, then the system (\ref{mainmgt0}) generates a strongly continuous semigroups on $\cH$. 
\end{theorem}
Moreover, introduce the parameter $\gamma=\alpha-{c^2\tau\over b}$, we have

\begin{theorem}[Kaltenbacher et.al. 2011, \cite{KLM11}]\label{KLM11_2}

Let $b>0, \alpha>0, c>0$. Let $$\hat{E}_1(t)=b||\Ahalf(u_t+{c^2\over b}u)||^2+\tau||u_{tt}+{c^2\over b}u_t||^2+{c^2\over b}\gamma ||u_t||^2,$$
$$\hat{E}_2(t)=\alpha||u_t||^2+c^2||\Ahalf u||^2,$$
and $$\hat{E}(t)=\hat{E}_1(t)+\hat{E}_2(t).$$

1. If $\gamma>0$, then the semigroup generated in Theorem \ref{KLM11_1} is exponentially stable on $\cH$. And there exist $\omega>0, C>0$ such that $$\hat{E}(t)\leq Ce^{-\omega t}\hat{E}(0), t>0.$$

2. If $\gamma=0$, the energy $\hat{E}_1(t)$ remains constant.
\end{theorem}

\begin{remark}
Theorems \ref{KLM11_1} and \ref{KLM11_2} were proved in \cite{KLM11}. Here we give a slightly different proof of Theorem \ref{KLM11_2}. In this process, the readers can familiarize themselves with the notations and the methodology that will be used in the later sections.
\end{remark}

\begin{proof} 

{\bf 1. Non-critical case: $\gamma>0$}.

We start with constructing energy functional and doing energy estimate.\\

Multiplying (\ref{mainmgt0}) by $u_{tt}$, we have

$${d\over dt}E_{01}(t)+R_{01}(t)=0$$
with
$$E_{01}(t)=[\tau ||u_{tt}||^2+b||\A ^{1/2}u_t||^2+2c^2(\A u, u_t)],$$
$$R_{01}(t)=2\alpha||u_{tt}||^2-2c^2||\A^{1/2}u_t||^2.$$

Multiply (\ref{mainmgt0}) by $u_t$ to get
$${d\over dt}E_{02}(t)+R_{02}(t)=0$$
with
$$E_{02}(t)=[c^2 ||\A^{1/2}u||^2+\alpha||u_t||^2+2\tau (u_{tt},u_t)],$$
$$R_{02}(t)=2b||\A^{1/2}u_t||^2-2\tau||u_{tt}||^2.$$
Sinces $\gamma=\alpha-{c^2\tau\over b}>0$, we can pick up a $k$ such that $${c^2\over b}<k<{\alpha\over\tau}.$$
Let $E_0(t)=E_{01}(t)+kE_{02}(t)$ be the {\bf natural energy functional}, we have
\begin{equation}\label{natural_energy}
\begin{split}
E_0(t)=b||\A^{1/2}u_t+{c^2\over b}\A^{1/2}u||^2+\tau||u_{tt}+ku_t||^2 \\
+k\tau({\alpha\over \tau}-k)||u_t||^2+c^2(k-{c^2\over b})||\A^{1/2}u||^2, 
\end{split}
\end{equation}
 and
\begin{eqnarray}\label{mgt_energy_identity}
&&{d\over dt}E_0(t)+R_0(t)=0, \\
R_0(t)&=&2\tau({\alpha\over\tau}-k)||u_{tt}||^2+2b(k-{c^2\over b})||\A^{1/2}u_t||^2.\nonumber
\end{eqnarray}

\begin{remark}\label{noncritical}
Note the coefficient of each term in $E_0$ and $R_0$ is strictly positive. This is possible only because of the existence of such a $k$ with ${c^2\over b}<k<{\alpha\over\tau}$. This is critical to generating an effective damper in the system, as we will see soon. This is also why we call it ``non-critical case'' in this work. For the ``critical case'', namely when $\gamma=\alpha-{c^2\tau\over b}=0$, no such a $k$ exists. In other words, the $k$ in the above calculations is forced to satisfied $k={c^2\over b}={\alpha\over\tau}$. So $R_0(t)=0$ and we have no damper in the energy equation, hence no energy decay can be obtained this way.
\end{remark}
\begin{remark}\label{natural_standard}
The {\bf natural energy functional} comes directly from the energy equality and normally looks messy due to complicated expressions or additional parameters. So we get it dressed up and introduce the equivalent {\bf standard energy functional}.
\end{remark}

We denote the {\bf standard energy functional} by 
\begin{equation}\label{standard_energy}
F_0(t)=||u_{tt}||^2+||\A^{1/2}u_t||^2+||\A^{1/2}u||^2.
\end{equation}
We have $E_0(t)\sim F_0(t)$, namely they are equivalent: there exists $C_1>0, C_2>0$ such that
$$C_1E_0(t)\leq F_0(t)\leq C_2E_0(t).$$
Indeed, recall $||u_t||\leq \lambda_0||\Ahalf u_t||$, this can be easily proved through the following lemma.
\begin{lemma}\label{square_lemma}
Let $C_0>0$. Then 
$$||f+g||^2+C_0||g||^2\sim ||f||^2+||g||^2.$$
\end{lemma}

\begin{proof}
Note $$||f+g||^2+C_0||g||^2={1\over {1+{C_0\over 2}}}||f||^2+2(f,g)+(1+{C_0\over 2})||g||^2$$
$$+(1-{1\over {1+{C_0\over 2}}})||f||^2+{C_0\over 2}||g||^2$$
$$\geq C_1(||f||^2+||g||^2), \mbox{  } C_1=\min\{(1-{1\over {1+{C_0\over 2}}}),{C_0\over 2}\}.$$ 
The other direction $||f+g||^2+C_0||g||^2\leq C_2(||f||^2+||g||^2)$ is trivial.
\end{proof}

Similarly, for $\hat{E}(t)$ in Theorem \ref{KLM11_2}, we can easily show $$\hat{E}(t)\sim F_0(t).$$ This simply means if we can show the exponential decay of $F_0(t)$, we get the same decay for $\hat{E}(t)$.

Generally, we want to prove the standard energy decays in certain way. Typically here we expect an exponential decay of $F_0$, namely, $F_0(t)\leq C_0e^{-Ct}$ for some positive constants $C_0, C$. One way of doing this is to construct Lyapunov function $L(t)$ based on the equation system, which is equivalent to energy $F_0(t)$; then one expects the Lyapunov function to show up in form of Gronwall inequality, which guarantee proper decay.

In the energy estimate, Gronwall inequalities take different forms. The most familiar one is in differential form, for non-negative function $L(t):[0,\infty)\rightarrow R^+$, $${d\over dt}L(t)\leq -CL(t)\Rightarrow L(t)\leq L(0)e^{-Ct}.$$ So if here we have 
$${d\over dt}L(t)\leq -C F_0(t),$$ or equivalently
$${d\over dt}L(t)\leq -C L(t),$$ then we get the exponential decay. 

In our case, it is natural to pick $L(t)=E_0(t)$. The problem is,  from (\ref{mgt_energy_identity}), we only have 
$${d\over dt}E_0(t)\leq -c(||u_{tt}||^2+||\A^{1/2}u_t||^2)\nleq -cE_0(t).$$
So the Gronwall inequality in differential form does not apply. We need it in more general form. The following lemma is taken from \cite[~p.1345]{ACS08}.

\begin{lemma}\label{gronwall_lemma}
Let $F(t)$ be a non-negative decreasing function on $[0,\infty)$ and if $C_1>0, T_0\geq 0$ exists such that $$\int_t^{\infty}F(s)ds\leq C_1F(t), t\geq T_0,$$
 then $$F(t)\leq F(0)\exp({1-{t\over C_1+T_0}}), \forall t\geq 0.$$
\end{lemma}


From the energy identity (\ref{mgt_energy_identity}), we know $E_0(t)$ is non-increasing. Moreover, integrate over $[t, T]$, we have
\begin{equation}\label{first_inequality} 
E_0(T)+C_3\int_t^T ||u_{tt}||^2+||\A^{1/2}u_t||^2ds\leq E_0(t),
\end{equation} for arbitrary $0\leq t\leq T$ and $C_3=\min\{2\tau({\alpha\over\tau}-k), 2b(k-{c^2\over b})\}$.
In particular we have
$$\int_t^T ||u_{tt}||^2+||\A^{1/2}u_t||^2ds\leq {1\over C_3}E_0(t).$$
Now if the following inequality holds for some constant $C_4>0$ independent of $t, T$, $$\int_t^T ||\A^{1/2}u||^2ds\leq C_4E_0(t),$$ we would arrive at 
$$\int_t^T E_0(s)ds\leq C_2\int_t^T F_0(s)ds\leq CE_0(t)$$ for arbitrary $0\leq t\leq T\leq \infty$. Hence take $T=\infty$, we would have the exponential decay of energy by Lemma \ref{gronwall_lemma}.

To show 
$$\int_t^T ||\A^{1/2}u||^2ds\leq C_4E_0(t),$$
multiply (\ref{mainmgt0}) by $u$,

$${d\over dt}{b\over 2}||\A^{1/2}u||^2+c^2||\A^{1/2}u||^2=\alpha||u_t||^2+{d\over dt}[{\tau\over 2}||u_t||^2-\tau(u_{tt},u)-\alpha(u_t,u)],$$ which gives
$${b\over 2}||\A^{1/2}u||^2|_t^T+c^2\int_t^T||\A^{1/2}u||^2ds=\alpha\int_t^T||u_t||^2ds+[{\tau\over 2}||u_t||^2-\tau(u_{tt},u)-\alpha(u_t,u)]|_t^T$$ 
$$\Rightarrow {b\over 2}||\A^{1/2}u(T)||^2+c^2\int_t^T||\A^{1/2}u||^2ds$$
$$\leq {b\over 2}||\A^{1/2}u(t)||^2+\alpha\int_t^T||u_t||^2ds+[{\tau\over 2}||u_t||^2-\tau(u_{tt},u)-\alpha(u_t,u)]|_t^T$$ 
$$\leq \alpha\int_t^T||u_t||^2ds+CF_0(t)+CF_0(T).$$ 
Since $||u_t||\leq \lambda_0 ||\A^{1/2}u_t||$, in view of (\ref{first_inequality}), the equivalence of $F_0(t)$ and $E_0(t)$ and the non-increasing property of $E_0(t)$, we have
$$\int_t^T||\A^{1/2}u||^2ds\leq C_4E_0(t).$$ So we get the exponential decay for $E_0(t)$,  hence for $F_0(t)$ and $\hat{E}(t)$, since they are all equivalent. The proof of first part of Theorem \ref{KLM11_2} with $\gamma>0$ is completed.

{\bf 2. Critical case: $\gamma=0$}.

In the critical case, $\gamma=\alpha-{c^2\tau\over b}=0$, all the calculations valid with 
$k=c^2/b=\alpha/\tau$, but the expressions shrink because of cancellations. We have
$$E_0^{cr}(t)=b||\A^{1/2}u_t+{c^2\over b}\A^{1/2}u||^2+\tau||u_{tt}+ku_t||^2,$$
$$=b||\A^{1/2}(u_t+ku)||^2+\tau||(u_t+ku)_t||^2,$$
 and

$${d\over dt}E_0^{cr}(t)+R_0^{cr}(t)=0,$$
$$R_0^{cr}(t)=2\tau({\alpha\over\tau}-k)||u_{tt}||^2+2b(k-{c^2\over b})||\A^{1/2}u_t||^2=0.$$
Here $E_0^{cr}(t)$ represents the energy functional in the critical case.  $R_0^{cr}(t)=0$ here suggests the nonexistence of damping. That is, we have $$E_0^{cr}(t)=E_0^{cr}(0).$$ It is not difficult to see that $E_0^{cr}(t)$ is exactly $\hat{E}_1(t)$ when $\gamma=0$. So the energy is conserved. This completes the proof of Theorem \ref{KLM11_2}. 
\end{proof}
\begin{remark}\label{substitution}
This is exactly the intuition provided by \cite{KLM11}: the MGT equation (\ref{mainmgt0}) can be written as $$\tau z_{tt}+b\A z+\gamma z_t={\gamma c^2\over b}u_t,$$ where $z(t)=u_t+{c^2\over b}u$. Since $\tau>0, b>0$, it is easy to see that, in the new variable $z$, (\ref{mainmgt0}) becomes wave equation when $\gamma=0$, hence no decay can be expected. This is exactly how MGT is connected to wave equation through the parameter $\gamma=\alpha-{c^2\tau\over b}$. On the other hand, when $\gamma>0$,  the left side of (\ref{mainmgt0}) becomes wave equation with frictional damping which intends to force exponential decay.
\end{remark}

\begin{remark}
From the spectrum viewpoint, this fact can be seen as well. Taking the prototype of (\ref{mainmgt0}) as an example where $\A=-\Delta$. After Fourier transform, letting $\hat{u}=\F[u]$, we have a third order ODE in time
$$\tau \hat{u}_{ttt}+\alpha \hat{u}_{tt}-b\xi^2\hat{u}_t-c^2\xi^2 \hat{u}=0.$$ It is not difficult to verify, the spectrum (roots of its characteristic equation $\tau r^3+\alpha r^2-b\xi^2 r-c^2\xi^2=0$) are strictly negative as long as $\gamma=\alpha-c^2\tau/b>0$; while one of the roots becomes zero when $\gamma=0$.
\end{remark}

From next section, we start to study MGT with memory terms.

\section{Memory of type 1: Proof of Theorem \ref{theorem1} }

We consider  the equation (\ref{mainmgt0}) with memory term of the first kind, so $w(t) = u(t) $. This leads to 

\begin{equation}\label{mainmgt1} 
\tau u_{ttt}+ \alpha u_{tt}+c^2\A u+b\A u_t -\int_0^tg(t-s)\A u(s)ds=0
\end{equation}

{\bf Proof of Theorem \ref{theorem1}}

%
%

Again we use multipliers and work through the energy estimate in several steps, where we define energy functional and show the decay of the solution. When multiplying (\ref{mainmgt1}) by $u_{tt}$ and $u_t$, only the memory terms need additional calculations because the other terms are the same as in last section. 

Note that here for the memory term, we have

a.) inner product with $u_{tt}$:

$$2(-g*\A u, u_{tt})={d\over dt}E_{11m}(t)+R_{11m}(t)$$
with 
\begin{equation}\label{e11m}
E_{11m}=-g'\circ\Ahalf u+g(t)||\Ahalf u||^2-2\int_0^tg(t-s)(\A u(s),u_t(t))ds,
\end{equation}
\begin{equation}\label{r11m}
R_{11m}=-g''\circ\Ahalf u+g'(t)||\Ahalf u||^2.
\end{equation}

b.) inner product with $u_t$:

$$2(-g*\A u, u_{t})={d\over dt}E_{12m}(t)+R_{12m}(t),$$ with
\begin{equation}\label{e12m}
E_{12m}=g\circ\Ahalf u-G(t)||\Ahalf u||^2.
\end{equation}
\begin{equation}\label{r12m}
R_{12m}=g'\circ\Ahalf u-g(t)||\Ahalf u||^2
\end{equation}
Now we are ready for the energy estimate.\\

{\bf Step} 1. Multiply by $u_{tt}$,  we have

$${d\over dt}E_{11}(t)+R_{11}(t)=0,$$ with
$$E_{11}(t)=E_{01}(t)+E_{11m}(t),$$
$$R_{11}(t)=R_{01}(t)+R_{1m}(t),$$
where $E_{11m}, R_{11m}$ are defined in (\ref{e11m}) and (\ref{r11m}).
Also recall $$E_{01}(t)=[\tau ||u_{tt}||^2+b||\A ^{1/2}u_t||^2+2c^2(\A u, u_t)]$$
and $$R_{01}(t)=2\alpha||u_{tt}||^2-2c^2||\A^{1/2}u_t||^2$$ as defined in the previous section.\\

{\bf Step} 2. Multiply by $u_t$, we have
$${d\over dt}E_{12}(t)+R_{12}(t)=0$$ with
$$E_{12}(t)=E_{02}(t)+E_{12m}(t),$$
$$R_{12}(t)=R_{02}(t)+R_{12m}(t).$$ Here $E_{12m}, R_{12m}$ are defined in (\ref{e12m}) and (\ref{r12m}) and
$$E_{02}(t)=[c^2 ||\A^{1/2}u||^2+\alpha||u_t||^2+2\tau (u_{tt},u_t)],$$
$$R_{02}(t)=2b||\A^{1/2}u_t||^2-2\tau||u_{tt}||^2.$$

{\bf Step} 3. Since $\gamma=\alpha-{c^2\tau\over b}>0$, we can pick a $k\in ({c^2\over b}, {\alpha\over\tau})$ and define the {\it natural energy} $E_1(t)=E_{11}(t)+kE_{12}(t)$. We have
$${d\over dt}E_1(t)+R_1(t)=0$$ with
$$E_1(t)=E_0(t)+E_{11m}+kE_{12m},$$
$$R_1(t)=R_0(t)+R_{11m}+kR_{12m}.$$
Here $E_0, R_0$ are from (\ref{natural_energy}) and (\ref{mgt_energy_identity}).\\

{\bf Step} 4. Now denote the {\it standard energy}

\begin{equation}\label{F_1}
F_1(t)=||u_{tt}||^2+||\Ahalf u_t||^2+||\Ahalf u||^2-g'\circ\Ahalf u.
\end{equation}
We claim that 

\begin{lemma}\label{first_energy_equiv}
$E_1(t)\sim F_1(t)$
\end{lemma}

\begin{proof}

From (\ref{natural_energy})-(\ref{mgt_energy_identity}), (\ref{e11m})-(\ref{r12m}), we have

\begin{eqnarray*}
E_1(t)=&&b||\Ahalf u_t||^2+2c^2(\A u, u_t)+c^2k||\Ahalf u||^2\\
&&+\tau||u_{tt}+ku_t||^2+k\tau({\alpha\over \tau}-k)||u_t||^2\\
&&-(g'-kg)\circ\Ahalf u+(g-kG(t))||\Ahalf u||^2\\
&&+2\int_0^tg(t-s)(\Ahalf u(t)- \Ahalf u(s), \Ahalf u_t(t))ds\\
&&-2G(t) (\A u(t), u_t(t))
\end{eqnarray*}
\begin{eqnarray*}
=&&b||\Ahalf u_t||^2+2(c^2-G(t))(\A u, u_t)+k(c^2-G(t))||\Ahalf u||^2\\
&&+\tau||u_{tt}+ku_t||^2+k\tau({\alpha\over \tau}-k)||u_t||^2\\
&&-(g'-kg)\circ\Ahalf u+g||\Ahalf u||^2\\
&&+2\int_0^tg(t-s)(\Ahalf u(t)- \Ahalf u(s), \Ahalf u_t(t))ds.
\end{eqnarray*}
\begin{eqnarray*}
=&&(b-{c^2-G(t)\over k})||\Ahalf u_t||^2\\
&&+{c^2-G(t)\over k}[||\Ahalf u_t||^2+2k(\A u, u_t)+k^2||\Ahalf u||^2\\
&&+\tau||u_{tt}+ku_t||^2+k\tau({\alpha\over \tau}-k)||u_t||^2\\
&&-(g'-kg)\circ\Ahalf u+g||\Ahalf u||^2\\
&&+2\int_0^tg(t-s)(\Ahalf u(t)- \Ahalf u(s), \Ahalf u_t(t))ds
\end{eqnarray*}
\begin{eqnarray}\label{E_1}
=&&{c^2-G(t)\over k}||\Ahalf u_t+k\Ahalf u||^2+(b-{c^2\over k})||\Ahalf u_t||^2\nonumber\\
&&+\tau||u_{tt}+ku_t||^2+k\tau({\alpha\over \tau}-k)||u_t||^2\nonumber\\
&&-g'\circ\Ahalf u+g||\Ahalf u||^2\\
&&+{G(t)\over k}||\Ahalf u_t||^2+kg\circ\Ahalf u\nonumber\\
&&+2\int_0^tg(t-s)(\Ahalf u(t)- \Ahalf u(s), \Ahalf u_t(t))ds.\nonumber
\end{eqnarray}
Since $$2|\int_0^tg(t-s)(\Ahalf u(t)- \Ahalf u(s), \Ahalf u_t(t))ds|\leq kg\circ\Ahalf u+{G(t)\over k}||\Ahalf u_t||^2,$$ the summation of the last three terms is nonnegative. So
\begin{eqnarray*}
E_1(t)\geq&&{c^2-G(t)\over k}||\Ahalf u_t+k\Ahalf u||^2+(b-{c^2\over k})||\Ahalf u_t||^2\\
&&+\tau||u_{tt}+ku_t||^2+k\tau({\alpha\over \tau}-k)||u_t||^2\\
&&-g'\circ\Ahalf u+g||\Ahalf u||^2.
\end{eqnarray*}
Recall $g>0, g'<0, G(t)=\int_0^tg(s)ds\leq\int_0^\infty g(s)ds=G(+\infty)< c^2$, by Lemma \ref{square_lemma}, we have $$E_1(t)\geq CF_1(t).$$
On the other hand, it is easy to see
\begin{eqnarray*}
E_1(t)=&&{c^2-G(t)\over k}||\Ahalf u_t+k\Ahalf u||^2+(b-{c^2\over k})||\Ahalf u_t||^2\\
&&+\tau||u_{tt}+ku_t||^2+k\tau({\alpha\over \tau}-k)||u_t||^2\\
&&-g'\circ\Ahalf u+g||\Ahalf u||^2\\
&&+{G(t)\over k}||\Ahalf u_t||^2+kg\circ\Ahalf u\\
&&+2\int_0^tg(t-s)(\Ahalf u(t)- \Ahalf u(s), \Ahalf u_t(t))ds\\
\leq&&C||\Ahalf u_t+k\Ahalf u||^2+C||\Ahalf u_t||^2\\
&&+C||u_{tt}+ku_t||^2+C||u_t||^2\\
&&-g'\circ\Ahalf u+C||\Ahalf u||^2\\
&&+C||\Ahalf u_t||^2+Cg\circ\Ahalf u\\
\leq && CF_1(t)
\end{eqnarray*} because of the boundedness of $g(t), G(t)$ and the Assumption \ref{assumption1}. Hence the equivalence of $E_1(t)$ and $F_1(t)$ is established.
\end{proof}

{\bf Step} 5.  Recall 
$${d\over dt}E_1(t)=-R_1(t),$$
\begin{eqnarray}\label{R_1}
R_1(t)=&&2\tau({\alpha\over\tau}-k)||u_{tt}||^2+2b(k-{c^2\over b})||\Ahalf u_t||^2\nonumber\\
&&+(g''-kg')\circ\Ahalf u+(kg-g')||\Ahalf u||^2\geq 0.
\end{eqnarray}

Since the right hand side does not contain a strictly negative multiple of $||\Ahalf u||^2$, we are one term away from the Gronwall inequality in derivative form. So we follow the idea in section 2, seeking the integral form.

From the above energy equality we have 
\begin{equation}
\begin{split}
C\int_t^T||u_{tt}||^2+||\Ahalf u_t||^2+(g''-kg')\circ\Ahalf u+(kg-g')||\Ahalf u||^2ds\\
+E_1(T)
\leq E_1(t)
\end{split}
\end{equation}
for any $0\leq t\leq T< \infty$. Note $g''\geq 0, g'\leq 0, g\geq 0$, in particular we have
\begin{equation}\label{third_inequality}
\int_t^T||u_{tt}||^2+||\Ahalf u_t||^2-g'\circ\Ahalf uds\leq CE_1(t).
\end{equation}

The Gronwall inequality in integral form
$\int_t^TE_1(s)ds\leq CE_1(t)$ would be established,  if only we can show $$\int_t^T||\Ahalf u||^2ds\leq CE_1(t).$$

{\bf Step 6}. To achieve the goal, we apply multiplier $u$ to (\ref{mainmgt})
$${d\over dt}{b\over 2}||\Ahalf u||^2+c^2||\Ahalf u||^2$$
$$=\alpha||u_t||^2+{d\over dt}[{\tau\over 2}||u_t||^2-\tau(u_{tt},u)-\alpha(u_t,u)]$$
$$+G(t)||\Ahalf u||^2-\int_0^tg(t-s)(\Ahalf (u(t)-u(s)),\Ahalf u(t))ds.$$ 
which gives
$$\int_t^T(c^2-G(t)-\epsilon G(t))||\Ahalf u||^2ds$$
$$\leq\alpha\int_t^T||u_t||^2ds+[{\tau\over 2}||u_t||^2-\tau(u_{tt},u)-\alpha(u_t,u)]|_t^T$$
$$+{1\over 4\epsilon}\int_t^Tg\circ\Ahalf uds-{b\over 2}||\Ahalf u||^2|_t^T.$$ 

Since $c^2>G(+\infty)\geq G(t)$, we can choose $\epsilon$ so small that $$(c^2-G(t)-\epsilon G(t))\geq C_{\delta}>0$$

Since 
$$\int_t^T||u_t||^2ds\leq \lambda_0^2\int_t^T||\Ahalf u_t||^2ds\leq CE_1(t),$$
$$||u_t||^2|_t^T\leq ||\Ahalf u_t(t)||^2+||\Ahalf u_t(T)||^2\leq C(E_1(t)+E_1(T))\leq CE_1(t),$$ 
$$|(u_{tt},u)||\leq |u_{tt}||^2+\leq |u_{t}||^2\leq CE_{1}(t),$$ 
$$|(u_{t},u)||\leq CE_{1}(t),$$ 
we have
$$C_{\delta}\int_t^T||\Ahalf u||^2ds$$
$$\leq CE_1(t)+C\int_t^Tg\circ\Ahalf uds\leq CE_1(t).$$ 
The last inequality is because of (\ref{third_inequality}) and the assumption $g'\leq -c_0g$, or equivalently $g\leq -{1\over c_0}g'$.

{\bf Step} 7. Now combining all the calculations above, we have $\int_t^TF_1(s)ds\leq CE_1(t)$, hence 
$$\int_t^TE_1(s)ds\leq C_1\int_t^TF_1(s)ds\leq CE_1(t).$$
Applying Lemma \ref{gronwall_lemma}, we have completed the proof of Theorem \ref{theorem1}. $\hfill\square$.

\begin{remark} Theorem \ref{theorem1} is for non-critical case. For the {\bf critical case}, the problem is open. Following all the calculations above, we arrive at (\ref{R_1}) but without terms in $R_1$ like $||u_{tt}||^2, ||\Ahalf u_t||^2$ as damper. In other words, only memory terms contribute to the damping. That creates challenging problems, which might require new multipliers or new methods.
\end{remark}

\section{Memory of type 2: Proof of Theorem \ref{theorem2} } 

We consider  the equation (\ref{mainmgt0}) with memory term of the second kind, so $w(t) = u_t(t) $. This leads to 

\begin{equation}\label{mainmgt2} 
\tau u_{ttt}+ \alpha u_{tt}+c^2\A u+b\A u_t -\int_0^tg(t-s)\A u_t(s)ds=0
\end{equation}

{\bf Proof of Theorem \ref{theorem2}} Next we work through the energy estimate.\\

{\bf Step} 1. Multiply by $u_{tt}$, we have

$${d\over dt}E_{21}(t)+R_{21}(t)=0$$
with 
$$E_{21}(t)=[\tau ||u_{tt}||^2+b||\A ^{1/2}u_t||^2+2c^2(\A u, u_t)]$$
$$+g\circ\Ahalf v-G(t)||\Ahalf v||^2$$
and 
$$R_{21}(t)=2\alpha||u_{tt}||^2-2c^2||\Ahalf u_t||^2$$
$$-g'\circ\Ahalf v+g(t)||\Ahalf v||^2.$$
Note $v(t)=u_t(t)$ here.

{\bf Step} 2. Multiply by $u_t$, we have
$${d\over dt}E_{22}(t)+R_{22}(t)=0$$ with
$$R_{22}(t)=2b||\Ahalf u_t||^2-2\tau||u_{tt}||^2-2G(t)||\Ahalf v||^2$$
$$+2\int_0^tg(t-s)(\A (v(t)-v(s)),u_t(t))ds$$ and
$$E_{22}(t)=[c^2 ||\Ahalf u||^2+\alpha||u_t||^2+2\tau (u_{tt},u_t)].$$

{\bf Step} 3. Since $\gamma=\alpha-{c^2\tau\over b}>0$, we can pick a $k$ such that ${c^2\over b}<k<{\alpha\over\tau}.$
Let $E_2(t)=E_{21}(t)+kE_{22}(t)$ be the natural energy, we have
$$E_2(t)=b||\Ahalf u_t||^2+2c^2(\A u, u_t)+c^2k||\Ahalf u||^2$$
$$+\tau||u_{tt}+ku_t||^2+k\tau({\alpha\over \tau}-k)||u_t||^2$$
$$+g\circ\Ahalf v-G(t)||\Ahalf v||^2.$$
$$=(b-G(t))||\Ahalf u_t||^2+2c^2(\A u, u_t)+c^2k||\Ahalf u||^2$$
$$+\tau||u_{tt}+ku_t||^2+k\tau({\alpha\over \tau}-k)||u_t||^2$$
$$+g\circ\Ahalf v.$$

 and
\begin{equation}\label{energy_2nd}
{d\over dt}E_2(t)+R_2(t)=0,
\end{equation}
$$R_2(t)=2\tau({\alpha\over\tau}-k)||u_{tt}||^2+2b(k-{c^2\over b})||\Ahalf u_t||^2$$
$$-g'\circ\Ahalf v+g(t)||\Ahalf v||^2$$
$$-2kG(t)||\Ahalf v||^2+2k\int_0^tg(t-s)(\A (v(t)-v(s)),u_t(t))ds.$$

{\bf Step} 4. Now denote $$F_2(t)=||u_{tt}||^2+||\Ahalf u_t||^2+||\Ahalf u||^2+g\circ\Ahalf v,$$
we claim that 
\begin{lemma}
$E_2(t)\sim F_2(t)$
\end{lemma}

\begin{proof}

Recall $G(t)=\int_0^tg(s)ds$ and $G(+\infty)\leq c_1<b-{c^2\over k}$.
\begin{eqnarray*}
E_2(t)
&=&(b-G(t))||\Ahalf u_t||^2+2c^2(\A u, u_t)+c^2k||\Ahalf u||^2\\
&+&\tau||u_{tt}+ku_t||^2+k\tau({\alpha\over \tau}-k)||u_t||^2\\
&+&g\circ\Ahalf v.
\end{eqnarray*}
By Assumption \ref{assumption2}, $G(t)\leq G(+\infty)< b-{c^2\over k}$, we have
$$(b-G(t))||\Ahalf u_t||^2+2c^2(\A u, u_t)+c^2k||\Ahalf u||^2\geq c(||\Ahalf u_t||^2+||\Ahalf u||^2)$$ for some positive constant $c$.
So we have $$E_2(t)\geq cF_2(t).$$
The other direction is trivial. So the equivalence of $E_2(t)$ and $F_2(t)$ is established.
\end{proof}

{\bf Step} 5.  Recall 
$$R_2(t)=2\tau({\alpha\over\tau}-k)||u_{tt}||^2+2b(k-{c^2\over b})||\Ahalf u_t||^2$$
$$-g'\circ\Ahalf v+g(t)||\Ahalf v||^2$$
$$-2kG(t)||\Ahalf v||^2+2k\int_0^tg(t-s)(\A (v(t)-v(s)),u_t(t))ds$$
$$=2\tau({\alpha\over\tau}-k)||u_{tt}||^2+2k(b-{c^2\over k}-G(t))||\Ahalf u_t||^2$$
$$-g'\circ\Ahalf v+g(t)||\Ahalf v||^2$$
$$+2k\int_0^tg(t-s)(\A (v(t)-v(s)),u_t(t))ds$$

Under the assumption $$g'\leq -c_0 g\Rightarrow -g'\geq c_0g $$
Now $$2k\int_0^tg(t-s)(\A (v(t)-v(s)),u_t(t))ds$$
$$\leq 2k[{\theta\over 2} G(t)||\Ahalf v||^2+{1\over 2\theta} g\circ \Ahalf v ]$$
$$=k\theta G(t)||\Ahalf v||^2+{k\over \theta} g\circ \Ahalf v $$

Under Assumption \ref{assumption2}, a straightforward calculation shows
$$R_2(t)\geq c(||u_{tt}||^2+||\Ahalf v||^2+g\circ\Ahalf v)$$ and
$$\int_t^T(||u_{tt}||^2+||\Ahalf v||^2+g\circ\Ahalf v)ds\leq cE_2(t).$$
The Gronwall inequality in integral form
$$\int_t^TE(s)ds\leq cE(t)$$ would be established,  if we can show $$\int_t^T||\Ahalf u||^2ds\leq cE(t).$$ 
And this is achieved by same idea as the previous section. We omit the calculation.
So the Theorem \ref{theorem2} is completed. $\hfill \square$.

%
%

\section{Memory of type 3: Proof of Theorems \ref{theorem31} and \ref{theorem32}  }

Now we are concerned with MGT with mixed memory. The equation reads

\begin{equation}\label{mgt31}
\tau u_{ttt}+ \alpha u_{tt}+c^2\A u+b\A u_t -\int_0^tg(t-s)\A w(s)ds=0,
\end{equation} where $w(s)=\lambda u(s)+u_t(s).$\\

{\bf Proof of Theorem \ref{theorem31}}\\

Let $E_3, R_3$ denote the energy and damper in this case. Recall how $\lambda$ is defined. We pick $k=\lambda$ and combine the calculations in the first two types of memory to get
$${d\over dt}E_3+R_3=0.$$
Here we have:
\begin{eqnarray}\label{e3r3}
E_3&=&({c^2\over k}-G(t))||\Ahalf w||^2+(b-{c^2\over k})||\Ahalf u_t||^2\nonumber\\
&&+\tau||u_{tt}+ku_t||^2+k\tau({\alpha\over \tau}-k)||u_t||^2\nonumber\\
&&+g\circ \Ahalf w,\nonumber\\
R_3&&=2(\alpha-k\tau)||u_{tt}||^2+2(bk-c^2)||\Ahalf u_t||^2\nonumber\\
&&+g||\Ahalf w||^2-g'\circ\Ahalf w.
\end{eqnarray}

Let \begin{equation}\label{F3}
F_3(t)=||\Ahalf u_t||^2+||\Ahalf u||^2+||u_{tt}||^2+g\circ \Ahalf w.
\end{equation}

Now we prove the Theorem \ref{theorem31}, namely, $F_3(t)$ decays exponentially.

\begin{proof}

Firstly, it is easy to see $E_3(t)\sim F_3(t)$. Also since $E_3'(t)=-R_3(t)\leq 0$, we have $E_3(t)$ non-increasing. Moreover,
\begin{eqnarray*}
&&C\int_t^T||u_{tt}||^2+||\Ahalf u_t||^2+g||\Ahalf w||^2+g\circ\Ahalf w(s)ds\\
&&\leq \int_t^T R_3(s)ds\leq E_3(t).
\end{eqnarray*}
In particular, 
$$\int_t^T||u_{tt}||^2+||\Ahalf u_t||^2+g\circ\Ahalf w(s)ds\leq CE_3(t).$$

Now if we could show 
\begin{equation}\label{aupart}
\int_t^T||\Ahalf u||^2ds \leq cE_3(t),
\end{equation}
we would complete the proof of Theorem \ref{theorem31} by Lemma \ref{gronwall_lemma}, following the same idea as before.

To show (\ref{aupart}), multiply (\ref{mgt31}) by $u$, we have 
\begin{eqnarray*}
0&=&\tau{d\over dt}(u_{tt},u)-\tau(u_{tt},u_t)\\
&&+\alpha{d\over dt}(u_t,u)-\alpha ||u_t||^2\\
&&+c^2||\Ahalf u||^2+{b\over 2}{d\over dt}||\Ahalf u||^2\\
&&-\int_0^tg(t-s)(\A w(s),u(t))ds\\
&=&\tau{d\over dt}(u_{tt},u)-{\tau\over 2}{d\over dt}||u_t||^2\\
&&+\alpha{d\over dt}(u_t,u)-\alpha ||u_t||^2\\
&&+c^2||\Ahalf u||^2+{b\over 2}{d\over dt}||\Ahalf u||^2\\
&&+\int_0^tg(t-s)(\A (w(t)-w(s)),u(t))ds\\
&&-\lambda G(t)||\Ahalf u||^2-{1\over 2}{d\over dt}[G(t)||\Ahalf u||^2]+{1\over 2}g(t)||\Ahalf u||^2.\\
\end{eqnarray*}
That leads to
\begin{eqnarray*}
&&(c^2-\lambda G(t)+{1\over 2}g(t))||\Ahalf u||^2\\
&\leq&-\tau{d\over dt}(u_{tt},u)+{\tau\over 2}{d\over dt}||u_t||^2\\
&&-\alpha{d\over dt}(u_t,u)+\alpha ||u_t||^2\\
&&-{b\over 2}{d\over dt}||\Ahalf u||^2\\
&&+\epsilon||\Ahalf u||^2+C_{\epsilon}g\circ\Ahalf w\\
&&+{d\over dt}{G(t)\over 2}||\Ahalf u||^2.
\end{eqnarray*}
Hence, 
\begin{eqnarray*}
&&\int_t^T(c^2-\lambda G(t)+{1\over 2}g(t)-\epsilon)||\Ahalf u||^2ds\\
&\leq&CE_3(t)+C\int_t^T (||u_t||^2+g\circ\Ahalf w(s))ds\\
&\leq&CE_3(t).
\end{eqnarray*}

Recall $G(+\infty)<{c^2\over \lambda}$, choose $\epsilon$ sufficiently small and $t$ sufficiently large so that $c^2-\lambda G(t)-g(t)-\epsilon$ is strictly positive, we get (\ref{aupart}), hence the exponential decay of $E_3(t)$. The same for $F_3(t)$ since they are equivalent.
\end{proof}

%
%

{\bf Proof of Theorem \ref{theorem32}}\\

We are considering $$\tau u_{ttt}+ \alpha u_{tt}+c^2\A u+b\A u_t -g*\A w=0, $$
with $\gamma=0, \lambda={c^2\over b}$. We shall prove $F_3^{cr}(t)  \triangleq ||\A^{1/2} w||^2 + ||w_t||^2 + g \circ  \A^{1/2} w$ decays exponentially.

\begin{proof}
Impose the conditions $k={c^2\over b}$ into (\ref{e3r3}), we get the energy estimate for the critical case 

$${d\over dt}E_3^{cr}+R_3^{cr}=0$$ with

\begin{eqnarray*}
E_3^{cr}=&&({c^2\over k}-G(t))||\Ahalf w||^2\\
&&+\tau||w_t||^2+g\circ \Ahalf w.\\
R_3^{cr}=&&g||\Ahalf w||^2-g'\circ\Ahalf w.
\end{eqnarray*}
It is easy to see $E_3^{cr}$ is equivalent to $F_3^{cr}$ since $G(+\infty)<{c^2\over k}$.

Because $\gamma=0$, as suggested in Remark \ref{substitution}, the equation can be written into
\begin{eqnarray}\label{zeq}
\tau w_{tt} + b \A w  -g*\A w= 0.
\end{eqnarray}

This is exactly the problem we dealt with in Lasiecka et.al.\cite{LW14}. From the analysis of the second order equations with memory, it is known that the system $(w,w_t) \in D(\A^{1/2} ) \times H $ decays exponentially. Adopting result, we have shown Theorem \ref{theorem32}.
\end{proof}


\begin{thebibliography}{1}
\bibitem{Adh00} S. Adhikari, {\it Structural Dynamic Analysis with Generalized Damping Models: Analysis}, Wiley-ISTE, (2013), 384 pp.
\bibitem{ACS08} F. Alabau-Boussouira, P. Cannarsa and D. Sforza, Decay estimates for second order evolution equations with memory. {\it Journal of Functional Analysis} {\bf 254}, (2008), 1342-1372.
\bibitem{AC09} F. Alabau-Boussouira and  P. Cannarsa, A general method for proving sharp energy decay rates for memory-dissipative evolution equations, {\it C. R. Acad. Sci. Paris, Ser.} I {\bf 347}, (2009), 867-872.
\bibitem{Ala10} F. Alabau-Boussouira, A unified approach via convexity for optimal energy decay rates of finite and infinite dimensional vibrating damped systems with applications to semi-discretized vibrating damped systems, {\it Journal of Differential Equations} {\bf 248}, (2010), 1473-1517.
\bibitem{Ala12} F. Alabau-Boussouira, On some recent advances on stabilization for hyperbolic equations. {\it Lecture Note in Mathematics 2048, CIME Foundation Subseries, Control of Partial Differential Equations, 1-100, Springer Verlag}, {\bf 2048}, (2012).

\bibitem{CCLW15} M. M. Cavalcanti, A D.D. Cavalcanti, I. Lasiecka and X. Wang, Existence and sharp decay rate estimates for a von Karman system with long memory, {\it Nonlinear Anal. Real World Appl.}, {\bf 22}, (2015), 289-306. 
\bibitem{CCM08} M.M. Cavalcanti, V.N.D. Cavalcanti and P. Martinez, General decay rate estimates for viscoelastic dissipative systems, {\it Nonlinear Analysis.}, {\bf 68} {1}, (2008), 177-193.

\bibitem{Che79} G. Chen, Control and stabilization for the wave equation in a bounded domain, {\it SIAM J. Control Optim.}, {\bf 17}, no. 1, (1979), 66-81.

\bibitem{CO03}M.M. Cavalcanti and H.P. Oquendo, Frictional versus viscoelastic damping in a semilinear wave equation, {\it SIAM J. Control Optim.}, {\bf 42} (2003), no.4, 1310-1324.

\bibitem{Cri79} D. G. Crighton, Model equations of nonlinear acoustics, {\it Ann. Rev. Fluid Mech.}, {\bf 11},
(1979), 11Ð33.

\bibitem{FP02} M. Fabrizio and  S. Polidoro, Asymptotic decay for some differential systems with fading memory, {\it Appl. Anal.}, {\bf 81} no. 6, (2002), 1245-1264.
\bibitem{HW10} X. Han and M. Wang, General decay rates of energy for the second order evolutions equations with memory, {\it Acta Applicanda, Math.}, {\bf 110}, (2010), 195-207.

\bibitem{HNR87}  W.J. Hrusa, J.A. Nohel and M. Renardy, {\it Mathematical Problems in Viscoelasticity}, Pitman Monographs and Surveys in Pure and Applied Mathematics {\bf 35}, Longman (1987), 273 pp. 
\bibitem{Jor09} P. M. Jordan, Nonlinear acoustic phenomena in viscous thermally relaxing fluids: Shock bifurcation and the emergence of diffusive solitons, {\it J. Acoust. Soc. Amer.}, {\bf 124}, 2491, (2008).

\bibitem{KLM11}  B. Kaltenbacher, I. Lasiecka and R. Marchand, well-posedness and exponential decay rates for the Moore-Gibson-Thompson equation arising in high intensity ultrasound, {\it Control Cybernet.}, {\bf 40}, (2011), no. 4, 971-988.
\bibitem{KLP12} B. Kaltenbacher, I. Lasiecka and M. Pospieszalska, Well-posedness and exponential decay of the energy in the nonlinear Jordan-Moore-Gibson-Thompson equation arising in high intensity ultrasound,
{\it Math. Models Methods Appl. Sci.}, {\bf 22}, (2012), no. 11.

\bibitem{Kor94} V. Komornik, {\it Exact Controllability and Stabilization. The Multiplier Method, Res. Appl. Math.}, Masson/Wiley, Paris/Chichester, (1994).

\bibitem{Kuz71} V. P. Kuznetsov, Equations of nonlinear acoustics, {\it Sov. Phys. Acoust.}, {\bf 16}, (1971), 467-470.

\bibitem{LMM13} I. Lasiecka, S. A. Messaoudi and M. I. Mustafa, Note on intrinsic decay rates for abstract wave equations with memory, {\it Journal of Mathematical Physics} {\bf 54}, 031504 (2013); doi: 10.1063/1.4793988.

\bibitem{LT93} I. Lasiecka  and D.  Tataru, Uniform boundary stabilization of semilinear wave equation with nonlinear boundary dissipation, {\it Differential and Integral Equations}, {\bf 6}, (1993), pp 507-533.

\bibitem{LW14} I. Lasiecka and X. Wang, Intrinsic Decay Rate Estimates for Semilinear Abstract Second Order Equations with Memory, {\it New Prospects in Direct, Inverse and Control Problems for Evolution Equations
Springer INdAM Series}, Vol 10, (2014), pp 271-303.

\bibitem{LW15} I. Lasiecka  and X. Wang, Moore-Gibson-Thompson equation with memory, part II: general decay of energy, submitted.

\bibitem{MMT12} R. Marchand, T. McDevitt and R. Triggiani, An abstract semigroup approach to the third-order Moore-Gibson-Thompson partial differential equation arising in high-intensity ultrasound: structural decomposition, spectral analysis, exponential stability, {\it Math. Methods Appl. Sci.}, {\bf 35}, (2012), no. 15, 1896-1929.

\bibitem{Mes08} S. A. Messaoudi,  General decay of solutions of a viscoelastic equation, {\it JMAA} {\bf 341} (2008), 1457-1467.

\bibitem{MM12} S. Messaoudi and M. Mustafa, General stability result for viscoelastic wave equations {\it Journal Mathematical Pgysics}, {\bf 53}, (2012).  
\bibitem{MG60} F. K. Moore and W. E. Gibson, Propagation of weak disturbances in a gas subject to relaxation effects, {\it J. Aerospace Sci. Tech.}, {\bf 27}, (1960), 117-127.
\bibitem{Ren00} M. Renardy, {\it Mathematical Analysis of Viscoelastic Flows}, CBMS-NSF Conference Series in Applied Mathematics {\bf 73}, SIAM (2000), 104 pp.
\bibitem{RS01} J. Rivera and A. Salvatierra, Asymptotic behaviour of the energy in partially viscoelastic materials, {\it Quart. Appl. Math.}, {\bf 59}, (2001), 557-578.

\end{thebibliography}
\end{document}